\documentclass[11pt,leqno]{amsart}
\usepackage{amssymb, amsmath}
\usepackage{pstricks,pst-node}

\renewcommand{\baselinestretch}{1.3}

\makeatletter
\renewcommand{\@seccntformat}[1]{\bf\@nameuse{the#1}.\quad}
\renewcommand\section{\@startsection{section}{1}%
            \z@{.7\linespacing\@plus\linespacing}{.5\linespacing}%
            {\normalfont\bfseries \boldmath}}
\renewcommand\subsection{\@startsection{subsection}{2}%
            \z@{.5\linespacing\@plus.7\linespacing}{-.5em}%
            {\normalfont\bfseries \boldmath}}
\renewcommand\subsubsection{\@startsection{subsubsection}{3}%
            \z@{.3\linespacing\@plus.5\linespacing}{-.5em}%
            {\normalfont\bfseries \boldmath}}
\makeatother

\theoremstyle{plain}
\newtheorem{theorem}{Theorem}

\newtheorem{lemma}[theorem]{Lemma}

\theoremstyle{definition}

\numberwithin{equation}{section}

\newcounter{listequation}

\def\note#1{{\small\tt <<#1>>}}  
\def\note#1{}              

\def\:{\colon}

\def\ZZ{\mathbb Z}

\def\RR{\mathbb R}
\def\CC{\mathbb C}

\def\Re{\operatorname{Re}}
\def\mod{\operatorname{mod}}

\def\:{\colon}

\def\JWw0S{{}^J W^{-w_0 S}}
\def\w0JWS{{}^{-w_0 J} W^S}
\def\fg{{\mathfrak g}}

\def\fp{{\mathfrak p}}

\def\fk{{\mathfrak k}}

\def\sl{{\mathfrak s\mathfrak l}}

\def\so{{\mathfrak s\mathfrak o}}

\def\SL{\operatorname{SL}}
\def\SO{\operatorname{SO}}
\def\Spin{\operatorname{Spin}}
\def\O{\operatorname{O}}

\def\isom{\cong}

\begin{document}
\title[The minimal representation of the conformal group]{The minimal Representation of the conformal group\\[5pt] and 
 classical 
solutions to the wave equation}

\author{Markus Hunziker}
\address{Department of Mathematics \\
            Baylor University\\ Waco, Texas  }
\email{\tt \{Markus\underline{\ }Hunziker, Mark\underline{\ }Sepanski, Ronald\underline{\ }Stanke\}@baylor.edu}

\author{Mark R. Sepanski}

\author{Ronald J. Stanke}

%




%


\maketitle

\begin{abstract}
Let $n$ be an integer $\geq 2$.  We consider the wave operator 
$\Box =-\partial_{t}^{2}+\sum_{i=1}^{n}\partial_{x_{i}}^{2}$,
where $(t,x)=(t,x_{1},\ldots,x_{n})$ are the canonical coordinates on 
Minkowski space $\mathbb{R}^{1,n}=\RR\times\RR^n$. 
Lie's prolongation algorithm
calculates the Lie algebra of infinitesimal symmetries of $\Box$ 
to be isomorphic to the conformal Lie algebra 
$\fg=\mathfrak{so}(2,n+1)$ plus an infinite dimensional piece
reflecting the fact that $\Box$ is linear. In particular,
$\ker \Box=\{f\in \mathcal{C}^{\infty}(\RR^{1,n})\mid \Box f=0\}$ 
is a representation of $\fg$. 
This Lie algebra action does not exponentiate to a global action of 
the conformal group $G=\SO(2,n+1)_0$ or any cover group.
However, for $n$ odd, it is known that $\ker \Box$ contains 
a nice $\fg$-invariant subspace that carries the minimal representation of $G$.
In this paper, we give a uniform realization of the 
minimal representation of a double cover of $G$ in $\ker \Box$ as a 
positive energy representation $\mathcal{H}^+$ for $n$ even and odd.   
Using this realization, we obtain an explicit orthonormal basis for $\mathcal{H}^+$
that is well behaved with respect to energy and angular momentum.
The lowest positive energy solution is, up to normalization, 
\begin{equation*}
   f(t,x)=\frac{1}{\sqrt{(1-it)^2+\Vert x\Vert^2}^{\,n-1}},
\end{equation*}
where $\sqrt{\ }$ denotes the principal branch of the square root.
Of special note,  for $n$ odd, all functions in our basis are rational functions.
Finally, using Fourier analysis with respect to this basis, we prove
that every classical real-valued solution to the wave equation
is the real part of a unique continuous element in the representation $\mathcal{H}^+$.
\end{abstract}

\section{Introduction}
\label{S:SL2}
\noindent
The minimal representations of the orthogonal groups are well known and much 
studied in representation theory.
Consequently, some of our methods and results overlap
with existing literature.  We mention a few of the most relevant 
references here. 
Important early work on the minimal representation of $\SO(4,4)$ was done by
B. Kostant starting with the paper \cite{Ko}.  B.~Binegar and R.~Zierau 
in \cite{BZ} then constructed the
minimal representation of $\SO(p,q)_{0}$ for $p+q$ even.   
Their model was based on the kernel of
the ultrahyperbolic wave operator $\Box_{p,q}$ acting on the space of
smooth functions on the cone
$
C^{p,q}=
\{(x,y)\in \mathbb{R}^{p,q}=\mathbb{R}^{p}\times \mathbb{R}%
^{q}\mid \Vert x\Vert =\Vert y\Vert \neq 0\}
$
of homogeneous degree $2-\frac{p+q}{2}$.  
In particular, for $p=2$ and $q=n+1$ with $n$ odd, their model was 
based on the kernel of the operator 
$\Box_{2,n+1}$ acting on homogenous function on the cone 
$C^{2,n+1}$ of degree $r=\frac{1-n}{2}$.
A more general study of homogenous
functions on generalized light cones was given by R. Howe and E. Tan in
\cite{HT} and connections to
dual pairs were studied by C.~Zhu and J.~Huang in \cite{ZH}. 
T.~Kobayashi and B.~Orsted made an exhaustive study of the minimal
representation of $\O(p,q)$ for $p+q$ even in \cite{KO1,KO2,KO3}.  
 They realized the representation as the kernel of the Yamabe operator 
$\Delta_{S^{p-1}\times S^{q-1}}$ acting on $\mathcal{C}^{\infty }(S^{p-1}\times
S^{q-1})$, as a certain subspace of the kernel of $\Box _{p-1,q-1}$ acting on $%
\mathcal{C}^{\infty }(\mathbb{R}^{p-1,q-1})$, and, via Fourier techniques, as $\mathcal{L}^{2}(%
{C}^{p-1,q-1})$. In particular, for $p=2$ and $q=n+1$ with $n$
odd, they realized the representation in a subspace of the kernel of the
usual wave operator $\Box=\Box _{1,n}$ acting on $\mathcal{C}^{\infty }(\mathbb{R}^{1,n})$. \
Finally, T.~Kobayashi and G.~Mano in \cite{KM} start with a
representation of $\SL(2,\mathbb{R})\times \SO(n)$ on $\mathcal{L}^{2}(\mathbb{R}^{n},%
\frac{dx}{\Vert x\Vert })$ and use a result of S.~Sahi \cite{Sa} to show the representation
extends to a double cover of $\SO(2,n+1)_{0}$ when $n$ is even.

This paper has three purposes.
Firstly, we give a uniform realization of the minimal representation of
a double cover $\widetilde{G}$  of $\SO(2,n+1)_0$ 
as a positive energy representation  $\mathcal{H}^+$ and a negative energy representation $\mathcal{H}^-$
in the kernel of the wave operator $\Box=\Box_{1,n}$ acting on  $\mathcal{C}^{\infty}(\RR^{1,n})$ for both $n$ even and odd.
Our  construction is related to a general construction due to T.~Enright and N.~Wallach
who in \cite{EW} showed for the Hermitian symmetric pairs of tube type 
every positive energy representation that is 
a reduction point  can be realized in the kernel of a generalized Dirac operator. 
In this context we are also able to  answer an open question concerning the irreducibility of the kernel.
Secondly,  we give an  explicit orthonormal basis of $\mathcal{H}^+$
relative to the Klein-Gordon inner product 
that is well behaved with respect to energy and angular momentum.
For $n$ odd, all functions in our basis are rational smooth solutions to the wave equation.
Thirdly, we use our explicit basis of $\mathcal{H}^+$ to study classical solutions to the wave equation.
We prove that every classical real-valued solution to the wave equation
is the real part of a unique smooth element in the representation $\mathcal{H}^+$.

For brevity, some proofs are only outlined or
omitted. Complete proofs of all our results will appear in a forthcoming 
paper.

\section{A Distinguished Copy of $\SL(2,\RR)$}
\label{S:SL2}
\noindent
Throughout the paper let $n$ be an integer $\geq 2$.
It is well known (c.f. \cite{Ol}) that Lie's prolongation algorithm
calculates the Lie algebra of infinitesimal symmetries 
of the wave operator 
$\Box =-\partial_{t}^{2}+\sum_{i=1}^{n}\partial_{x_{i}}^{2}$
on $\RR^{1,n}=\RR\times\RR^n$
to be isomorphic to the conformal Lie algebra 
$\fg=\mathfrak{so}(2,n+1)$ plus an infinite dimensional piece
reflecting the fact that $\Box$ is linear.
The Lie algebra $\fg$  contains a distinguished copy of $\sl(2,\RR)$
spanned by the triple
\begin{equation}
\begin{aligned}
  h\ \ &=\textstyle{2\left(r-t\partial_t-\sum_i x_i\partial_{x_i}\right)},\\
  e^+&=-\partial_t, \\
  e^-&=\textstyle{2t\left(r-t\partial_t-\sum_i x_i\partial_{x_i}\right)+q(t,x)\partial_t},
  \end{aligned}
\end{equation}
where $r=\frac{1-n}{2}$ and $q(t,x)=-t^2+\Vert x\Vert^2$. 
The centralizer of this copy of  $\sl(2,\RR)$ in $\fg$ is the Lie algebra $\so(n)$ 
spanned by the infinitesimal rotations 
$-x_j\partial_{x_i}+x_i\partial_{x_j}$, $1\leq i<j\leq n$.
If we define $\so(2,n+1)$ with respect to the symmetric bilinear form
given by $\operatorname{diag}(-1,-1,1,\ldots,1)$ then 
$\sl(2,\RR)\times \so(n)$ embeds 
block diagonally as the Lie subalgebra $\so(2,1)\times \so(n)$. 
Corresponding to this Lie subalgebra we have
the block diagonally embedded subgroup $\SO(2,1)_0\times \SO(n)$
of the conformal group $G=\SO(2,n+1)_0$. 
Also diagonally embedded in $G$ is the maximal compact subgroup
$K=\SO(2)\times \SO(n+1)$.
The group $K$ has a double cover $\widetilde{K}\isom \SO(2)\times \SO(n+1)$   
with covering map $\pi : \widetilde{K} \rightarrow K$ given by
$(R_{\frac{\varphi}{2}},k) \mapsto (R_{\varphi},k)$, where
$R_\alpha \in \SO(2)$ is the rotation by  $\alpha$ and $k \in \SO(n+1)$.
It is then possible to extend $\widetilde{K}$ to a connected Lie group $\widetilde{G}$
with Lie algebra $\mathfrak{g}$ and extend $\pi $ to a map
 $\pi :\widetilde{G}\rightarrow G$ so that $\widetilde{G}$ is a double cover of $G$ with
covering map $\pi $. The double cover $\widetilde{G}$ contains 
a subgroup isomorphic to $\SL(2,\RR)\times \SO(n)$ such that the restriction of the 
covering map  $\pi : \widetilde{G}\rightarrow G$ to the 
first factor is the  map $\SL(2,\RR)\isom \Spin(2,1)_0\rightarrow \SO(2,1)_0$.
The group $\widetilde{G}$ and its distinguished  subgroup
$\SL(2,\RR)\times \SO(n)$ will play a pivotal role in our
study of solutions to the wave equation. 
The starting point is the formula
\begin{equation}\label{Omega}
\Omega_{\SL(2)}-\Omega_{SO(n)}-r(r+1) =\Vert x\Vert^2 \Box,
\end{equation}
where $\Omega_{\SL(2)}=\frac{1}{4}h^2+\frac{1}{2}(e^+e^- + e^-e^+)$ and
$\Omega_{\SO(n)}=- \sum_{i<j}(-x_j\partial_{x_i}+x_i\partial_{x_j})^2$
are the Casimir elements of the Lie subalgebras $\sl(2,\RR)$ and
$\so(n)$ of $\frak{g}$ with respect to the bilinear form 
$B(X,Y)=\frac{1}{2}\operatorname{tr}(XY)$ on $\fg$.

\section{Degenerate Principal Series Representations}
\label{S:DegPrSeries}
\noindent
By construction via the prolongation algorithm, the Lie algebra $\fg$ acts naturally on
$\ker \Box=\{f\in \mathcal{C}^{\infty}(\RR^{1,n})\mid \Box f=0\}$.
This Lie algebra action does not exponentiate to an action of the group $G$ or its
double cover $\widetilde{G}$. We will now show
how to obtain a group action on a subspace of $\ker \Box$ given by restricting smooth sections 
in the kernel of the differential operator $\Omega=\Omega_{\SL(2)}-\Omega_{\SO(n)}-r(r+1)$,
$r=\frac{1- n}{2}$, acting on an equivariant line bundle on  the conformal compactification  of $\RR^{1,n}$.

The conformal compactification of $\RR^{1,n}$ arises naturally in our setting 
as follows. The eigenvalues of $\operatorname{ad}(h)$ on $\mathfrak{g}$ are $\{2,0,-2\}$ which
gives rise to a pair of opposite maximal parabolic subgroups $Q^{\pm }$ of $G=\SO(2,n+1)_0$ 
with Langlands decompositions $Q^{\pm }=MAN^{\pm }$ for which $M\cong \SO(1,n)$,
$A\cong \mathbb{R}_{>0}$, and $N^{\pm }\cong \mathbb{R}^{1,n}$.
Similarly, we have a pair  of opposite maximal parabolic subgroups $\widetilde{Q}^{\pm }$ of $\widetilde{G}$ 
with Langlands decompositions
 $\widetilde{Q}^{\pm }=\widetilde{M}\widetilde{A}\widetilde{N}^{\pm }$ 
 for which $\widetilde{M}_{0}\cong \SO(1,n)_{0}$, $\widetilde{A}\cong \mathbb{R}_{>0}$, and
 $\widetilde{N}^{\pm }\cong \mathbb{R}^{1,n}$. The group $\widetilde{M}$ has four
connected components with $\widetilde{M}/\widetilde{M}_{0}\cong \mathbb{Z}_{4}$.
Identifying $\RR^{1,n}$ with $N^+\isom \widetilde{N}^+$, we obtain an embedding
$\RR^{1,n} \hookrightarrow G/Q^-\isom \widetilde{G}/\widetilde{Q}^-$ with open dense image.
This is the conformal compactification of $\RR^{1,n}$. 

Turning to line bundles, recall that
the $\widetilde{G}$-equivariant line bundles on $\widetilde{G}/\widetilde{Q}^-$
are of the form
$\mathcal{L}_\chi =\widetilde{G}\times_{\widetilde{Q}^-} \CC_\chi$, where $\chi$ is a character
of $\widetilde{Q}^-$. The space of global sections of $\mathcal{L}_\chi$
is the {\it degenerate principal series representation} $\operatorname{Ind}_{\widetilde{Q}^-}^{\widetilde{G}}(\chi)$
defined by 
\begin{equation*}
\operatorname{Ind}_{\widetilde{Q}^-}^{\widetilde{G}}(\chi)=
\left\{ \phi\in \mathcal{C}^{\infty }(\widetilde{G})
\mid\phi(gq^{-})=
\chi ^{-1}
(q^{-})\phi(g)\quad \forall g\in \widetilde{G}, q^-\in \widetilde{Q}^{-}
\right\}
\end{equation*}
with $\widetilde{G}$-action given by 
$( g\cdot \phi)(x)=\phi(g^{-1}x)$.
Restricting functions to $\widetilde{N}^+\isom \RR^{1,n}$,
we obtain the so-called {\it noncompact picture\/} $\mathcal{I}'_{\chi}\subset \mathcal{C}^{\infty}(\RR^{1,n})$
of $\operatorname{Ind}_{\widetilde{Q}^-}^{\widetilde{G}}(\chi)$.

A character  $\chi$ of $Q^-$  is determined
by a discrete parameter $m\in \ZZ_{4}$ and a continuous parameter $r\in \CC$
as follows.
Let $\gamma _{0}:\mathbb{Z}_{4}\rightarrow S^{1}$ be given by $\gamma_{0}(j)=i^{j}$. 
Identifying $\widetilde{M}/\widetilde{M}_{0}$ with $\mathbb{Z}_{4}$, 
let $\gamma :\widetilde{M}\rightarrow S^{1}$ be the pull back of $%
\gamma _{0}$ given by mapping $\widetilde{M}$ to its component group and
then applying $\gamma _{0}$. For $m\in \mathbb{Z}_{4}$ and $r\in \mathbb{C}$, 
we then define a character
$\chi_{m ,r}:\widetilde{Q}^{-}\rightarrow \mathbb{C}^{\times }$ 
by 
$\chi_{m ,r}(q^-) =\gamma ^{m}\big(q^-_{\widetilde{M}}\big)e^{rs}$,
where $q^-=q^-_{\widetilde{M}}\,q^-_{\widetilde{A}}\,q^-_{\widetilde{N}^-}$
is written with respect to the Langlands decomposition with 
$q^-_{\widetilde{A}}=\exp_{\widetilde{G}}(s h)$, $s\in \RR$.
We remark that every character of $\widetilde{Q}^{-}$ is of the form $\chi_{m,r}$
and $\operatorname{Ind}\nolimits_{\widetilde{Q}^{-}}^{\widetilde{G}}( \chi_{m ,r}) $ 
descends to a representation of $G$ if and only if $m$ is even.
Writing $\mathcal{I}'_{m,r}\subset \mathcal{C}^{\infty}(\RR^{1,n})$
for the noncompact picture of $\operatorname{Ind}\nolimits_{\widetilde{Q}^{-}}^{\widetilde{G}}( \chi_{m ,r})$,
the Lie algebra $\fg$ acts on $\mathcal{I}'_{m,r}$ by first order differential operators. 
These differential operators only depend on the parameter $r$
and hence can be calculated explicitly by choosing $m=0$ and viewing $\mathcal{I}'_{m,r}$ as a representation of $G$.
It turns out that for the special parameter $r=\frac{1-n}{2}$, the $\fg$-action on $\mathcal{I}'_{m,r}$ is given by exactly the
same differential operators as the $\fg$-action on $\mathcal{C}^\infty(\RR^{1,n})$
that is obtained via the prolongation algorithm. In particular,
the subspace $\ker \Box \subset \mathcal{I}_{m,r}'$ is  $\fg$-invariant and hence 
$\widetilde{G}$-invariant. (Here we slightly abuse  language and write $\ker \Box \subset \mathcal{I}_{m,r}'$
for $\ker \Box\vert_{\mathcal{I}_{m,r}'}$.)
In light of formula (\ref{Omega}),  $\ker \Box\subset \mathcal{I}'_{m,r}$  corresponds to $\ker \Omega \subset 
 \operatorname{Ind}\nolimits_{\widetilde{Q}^{-}}^{\widetilde{G}}( \chi_{m ,r})$
under the isomorphism $\mathcal{I}'_{m,r} \isom   
\operatorname{Ind}\nolimits_{\widetilde{Q}^{-}}^{\widetilde{G}}( \chi_{m ,r})$ for $r=\frac{1-n}{2}$.
Before we study the $\widetilde{G}$-representations $\ker \Omega \subset 
 \operatorname{Ind}\nolimits_{\widetilde{Q}^{-}}^{\widetilde{G}}( \chi_{m ,r})$
in  detail, we will give three  concrete realizations of $\operatorname{Ind}_{\widetilde{Q}^-}^{\widetilde{G}}(\chi_{m,r})$
in the following sections.

\section{Geometric Picture}

\noindent
Define the cone $C^{2,n+1}=\{(a,b)\in \mathbb{R}^{2,n+1}\mid
\left\Vert a\right\Vert =\left\Vert b\right\Vert \neq 0\}$. $\ $The group $G$
acts on $C^{2,n+1}$ by matrix multiplication on the left. Now
consider  the double cover of $C^{2,n+1}$ defined by $\widetilde{%
C}^{2,n+1}=C^{2,n+1}$ with covering map $\eta:\widetilde{%
C}^{2,n+1}\rightarrow C^{2,n+1}$ given by $\eta(
\lambda \sin \frac{\varphi }{2},\lambda \cos \frac{\varphi }{2},\lambda
b) =( \lambda \sin \varphi ,\lambda \cos \varphi ,\lambda
b) $ where $\lambda \in \RR_{>0}$, $\varphi \in \mathbb{R}$, and 
$b\in S^{n}$.
A priori, $\widetilde{C}^{2,n+1}$ carries only an action of 
$\widetilde{K}=\SO(2)\times \SO(n+1)$ given by matrix multiplication on the left that is compatible with $%
\eta$ in the sense that $\eta(k\cdot c) =\pi (k)\cdot \eta(c)$ for $k\in \widetilde{K}$ and $c\in \widetilde{%
C}^{2,n+1}$. However, using the diffeomorphisms%
\begin{equation*}
\widetilde{G}/( \widetilde{M}_{0}\widetilde{N}^{-}) \overset{%
\cong }{\longleftarrow }\widetilde{K}/( \widetilde{K}\cap \widetilde{M}%
_{0}) \times \widetilde{A}\overset{\cong }{\longrightarrow }\widetilde{%
C}^{2,n+1}\text{,}
\end{equation*}%
it follows that the $\widetilde{K}$-action   on $\widetilde{C}%
^{2,n+1}$ extends to a $\widetilde{G}$-action. (This action is {\it not\/} the restriction of a linear action.) Furthermore,
there is a commutative diagram%
\begin{equation*}
\begin{array}{ccc}
\widetilde{G} & \circlearrowright & \widetilde{C}^{2,n+1} \\ 
\begin{array}{cc}
\pi & \downarrow%
\end{array}
&  & 
\begin{array}{cc}
\downarrow & \eta%
\end{array}
\\ 
G & \circlearrowright & C^{2,n+1}\text{,}%
\end{array}%
\end{equation*}%
i.e., $\eta( g\cdot c) =\pi (g)\cdot \eta(c)$ for $g\in \widetilde{G}$, 
$c\in \widetilde{C}^{2,n+1}$. As a  side remark, we note that it is possible to use
this diagram to realize $\widetilde{G}$ as the elements $\widetilde{g}\in 
\operatorname{Diff}\nolimits^{\infty }(\widetilde{C}^{2,n+1})$ so that
there exists a $g\in G$ satisfying $\eta( \widetilde{g}(c)) =g\cdot \eta(c)$ for all $c\in \widetilde{C}^{2,n+1}$.

Now let $w$ be the block diagonally embedded element of $\SO(2)\times \O(n+1)$
given by $(R_{\frac{\pi }{2}},-I_{n+1})$ and acting on $\widetilde{{C%
}}^{2,n+1}$ by matrix multiplication on the left. 
The action of $w$ commutes with the $\widetilde{G}$-action and  
there is a $\widetilde{G}$--equivariant diffeomorphism $\widetilde{G}/( 
\widetilde{M}\widetilde{N}^{-}) \overset{\cong }{\longrightarrow }%
\widetilde{C}^{2,n+1}/\left\langle w\right\rangle $. From this we
obtain a more geometric and isomorphic realization of 
$\operatorname{Ind}\nolimits_{\widetilde{Q}^{-}}^{\widetilde{G}}( 
\chi_{m ,r}) $ as%
\begin{equation*}
\mathcal{I}_{m ,r}=\left\{ \mathcal{\phi }\in \mathcal{C}^{\infty }(\widetilde{%
C}^{2,n+1})\mid \mathcal{\phi }(w\cdot c)=i^{-m}\mathcal{\phi }(c)\text{
and }\mathcal{\phi }(\lambda c)=\lambda ^{r}\mathcal{\phi }(c)\text{, }%
\lambda \in \mathbb{R}_{>0}\right\}
\end{equation*}%
with $\widetilde{G}$ action given by $( g\cdot \mathcal{\phi })
(c)=\mathcal{\phi }(g^{-1}\cdot c)$, $g\in \widetilde{G}$, $c\in \widetilde{%
C}^{2,n+1}$. 

\section{Noncompact Picture}
\noindent
By letting $\widetilde{N}^{+}\isom \RR^{1,n}$ act on the base point $(0,1,-1,0,\ldots
,0)\in \widetilde{C}^{2,n+1}$, we obtain a map $\iota :\mathbb{R}%
^{1,n}\rightarrow \widetilde{C}^{2,n+1}$. Explicitly, 
$\iota(t,x)=(2t,1+q(t,x),-1+q(t,x),2x)$,
where $q(t,x)=-t^2+\Vert x\Vert^2$ as before.
We then define the noncompact picture
of $\mathcal{I}_{m ,r}$ as
\begin{equation*}
\mathcal{I}_{m ,r}^{\prime }=\left\{ f\in \mathcal{C}^{\infty }(\mathbb{R}^{1,n})\mid
\exists \mathcal{\phi }\in \mathcal{I}_{m ,r}\text{ so }f=\mathcal{\phi }%
\circ \iota \right\} \text{.}
\end{equation*}
This definition of $\mathcal{I}_{m ,r}^{\prime }$ is equivalent to the definition that was
given earlier.


\begin{lemma}\label{Lemma: j}
Let $c=( \lambda \sin \frac{\varphi }{2},\lambda \cos \frac{\varphi }{2},\lambda b)\in  \widetilde{C}^{2,n+1}$ with 
$\lambda \in \RR_{>0}$, $\varphi\in \RR$, and
$b=(b_{0},b_{1},\ldots ,b_{n})\in S^{n}$. If  $\phi \in \mathcal{I}_{m ,r}$ and $f=\phi \circ \iota$, then
\begin{equation*}
\phi (c)=i^{m j}\left( \frac{\lambda \left\vert \cos \varphi
-b_{0}\right\vert }{2}\right) ^{r}
f\left( \frac{\sin \varphi }{\cos \varphi
-b_{0}},\frac{b_{1}}{\cos \varphi -b_{0}},\cdots ,\frac{b_{n}}{\cos \varphi
-b_{0}}\right) ,  \label{eqn:  ref back for j}
\end{equation*}
where $j=j(\varphi,b_0)\in \ZZ_4$ is given by
 \begin{equation*}
  j=
  \begin{cases}
  0, &\mbox{if\quad $\cos\varphi-b_{0}>0$ and $\frac{\varphi}{2}\in(-\frac{\pi}{2},\frac{\pi}{2}) \mod 2\pi$ },\\
  1, &\mbox{if\quad $\cos\varphi-b_{0}<0$ and $\frac{\varphi}{2}\in(0,\pi) \mod 2\pi$},\\
  2, &\mbox{if\quad $\cos\varphi-b_{0}>0$ and $\frac{\varphi}{2}\in(\frac{\pi}{2},\frac{3\pi}{2}) \mod 2\pi$},\\
  3, &\mbox{if\quad $\cos\varphi-b_{0}<0$ and $\frac{\varphi}{2}\in(\pi,2\pi)\mod 2\pi$}.\\
  \end{cases}
\end{equation*}
\end{lemma}

\begin{proof}
The formula follows  from the definitions by a direct calculation. 
\end{proof}

Calculating the explicit action of $\widetilde{G}$ on $\mathcal{I}_{m ,r}^{\prime }$ 
is rather subtle, but always possible by using the lemma. As an
important example, we give (without proof) the action of the subgroup $\SL(2,\mathbb{R})\times \SO(n) 
\subset \widetilde{G}$.
\begin{theorem}
Let $g=(
\left(\begin{smallmatrix}
a & b \\ 
c & d%
\end{smallmatrix}\right)
 ,k)\in \SL(2,\mathbb{R})\times \SO(n)$ and $f\in \mathcal{I}%
_{m ,r}^{\prime }$. Then%
\begin{equation*}
( g\cdot f) (t,x)=( \sqrt{\operatorname{sgn}\delta }\,)
^{m}\left\vert \delta \right\vert ^{r}f\left( \frac{(-b+dt)(a-ct)+cd\left%
\Vert x\right\Vert ^{2}}{\delta },\frac{xk}{\delta }\right)
\end{equation*}%
where $\delta =(a-ct)^{2}-c^{2}\left\Vert x\right\Vert ^{2}$ and $\sqrt{%
\operatorname{sgn}\delta }$ is defined as $+1$ if $\delta >0$ and $%
a-ct>c\left\Vert x\right\Vert $, as $-1$ if $\delta >0$ and $%
a-ct<c\left\Vert x\right\Vert $, and as $i$ if $\delta <0$. (See also Fig.~1.)
\qed
\end{theorem}

%

\begin{figure}[ht]\label{F:light cone}
\centering
\begin{pspicture}(0,-2.5)(0,7.5)
\psset{linewidth=.5pt,labelsep=8pt,nodesep=0pt}
\psline{->}(3,3.5)(2,3)
\psellipse[border=2pt](0,5.05)(2,.5)
\psarc[border=2pt]{<-}(-2,4.8){1}{10}{70}
\psarc[border=2pt]{->}(-2,1.2){1}{-70}{-10}
\psline{->}(0,0)(0,6)
\psellipse[border=2pt](0,0.95)(2,.5)
\psline[border=2pt](0,1)(0,3)
\psline[border=2pt](0,5)(0,5.8)
\psline{->}(0,0)(-1,-2)
\psline{->}(0,0)(3,0)
\psline(0,3)(1.95,4.95)
\psline(0,3)(-1.95,4.95)
\psline(0,3)(-1.95,1.05)
\psline(0,3)(1.95,1.05)
\uput[l](-1.8,5.8){$\sqrt{\operatorname{sgn}(\delta)}=+1$}
\uput[l](-1.8,0.2){$\sqrt{\operatorname{sgn}(\delta)}=-1$}
\uput[r](3,3.7){$\sqrt{\operatorname{sgn}(\delta)}=i$}
\uput[r](0,3){$\frac{a}{c}$}
\uput[u](0,6){$t$}
\uput[r](3,0){$x_{2}$}
\uput[d](-1.1,-2){$x_{1}$}
\end{pspicture}
\caption{The values of $\sqrt{\operatorname{sgn}(\delta)}$}
\end{figure}




\section{Compact Picture}

\noindent
The final realization of $\operatorname{Ind}\nolimits_{\widetilde{Q}^{-}}^{%
\widetilde{G}}( \chi_{m ,r}) $ is essentially given by
restricting the elements of $\mathcal{I}_{m ,r}$ to $S^{1}\times
S^{n}\subseteq \widetilde{C}^{2,n+1}$, i.e., the so-called {\it compact
picture\/}. Because of the prominent role of the group $\SL(2,\RR)\times \SO(n)$, we use
spherical coordinates on $S^{n}$ and pull $S^{1}$ back to $\mathbb{R}$ by
the map $\sigma :\mathbb{R}\times \mathbb{R}\times S^{n-1}\rightarrow
S^{1}\times S^{n}$ given by $\sigma (\varphi ,\theta ,\widehat{x})=( \sin 
\frac{\varphi }{2},\cos \frac{\varphi }{2},-\cos \theta ,\widehat{x}\sin
\theta ) $. 
We then define
$$
\mathcal{I}_{m ,r}^{\prime \prime
}=\{F\in \mathcal{C}^{\infty }(\mathbb{R}\times \mathbb{R}\times S^{n-1})\mid
\mbox{$F=\phi \circ \sigma $, some $\phi \in \mathcal{I}_{m ,r}$}\}
$$ and give it a $%
\widetilde{G}$-structure so that the map $\phi \rightarrow \phi \circ \sigma $
is an isomorphism. 
If $\phi \in \mathcal{I}_{m ,r}$ corresponds to $f\in \mathcal{I}_{m ,r}^{\prime }$ 
and to $F\in \mathcal{I}_{m ,r}^{\prime \prime }$, then it follows from the definitions
and Lemma~\ref{Lemma: j} that
\begin{equation}\label{noncompact to compact}
F(\varphi ,\theta ,\widehat{x}) =i^{mj}\left\vert \frac{\cos \varphi +\cos
\theta }{2}\right\vert ^{r}f\left( \frac{\sin \varphi }{\cos \varphi +\cos
\theta },\frac{\widehat{x}\sin \theta }{\cos \varphi +\cos \theta }\right)\ \mbox{and}\qquad\quad
\end{equation}
\begin{equation}\label{compact to noncompact}
f(t,x) =\lambda (t,x)^{r}\,F\left( \operatorname{sgn}(t)\cos ^{-1}\left( \frac{%
1+q(t,x)}{\lambda (t,x)}\right) ,\cos ^{-1}\left( \frac{1-q(t,x)%
}{\lambda (t,x)}\right) ,\frac{x}{\left\Vert x\right\Vert }\right),
\end{equation}
where $j$ is  as in Lemma~~\ref{Lemma: j}  with $b_0=-\cos\theta$
and $\lambda(t,x)=\left((1-q(t,x))^{2}+4\Vert x\Vert^{2}\right)^{\frac{1}{2}}$.

\begin{theorem}\label{Omega compact}
For $r=\frac{1-n}{2}$,  we have the following identity of differential  operators on $\mathcal{I}_{m,r}''$: 
\begin{equation*}
   \Omega_{\SL(2)}-\Omega_{\SO(n)}-r(r+1)= \sin^2\theta\, \left(\Omega_{\SO(2)}-\Omega_{\SO(n+1)}-r^2\right).
\end{equation*}
\end{theorem}

\begin{proof}
Using the definitions and (\ref{compact to noncompact})
it is easy to show that $\Omega_{\SL(2)}$ acts on $\mathcal{I}''_{m,r}$ by the formula
$
  \Omega_{\SL(2)}= r(1+r) -r^{2}\sin^{2}\theta  -2r\cos\theta\sin\theta\, \partial_{\theta}+\sin^{2}\theta
\left(-\partial_{\varphi}^{2}+\partial_{\theta}^{2}\right)$.
Moreover, $\Omega_{\SO(2)}=-\partial_\varphi^2$, 
$\Omega_{\SO(n+1)}=-\Delta_{S^n}$, and $\Omega_{\SO(n)}=-\Delta_{S^{n-1}}$.
The theorem then follows 
from the well known recursive formula
\begin{equation}\label{spherical laplacian}
\Delta_{S^n} =  \partial_{\theta}^{2} +(n-1) \cot \theta\, \partial_{\theta} -
\csc^2\theta\, \Delta_{S^{n-1}}
%
\end{equation}
for the spherical Laplacian in spherical coordinates.
\end{proof}

\section{$\widetilde{K}$-Types}

\noindent 
For $r=\frac{1-n}{2}$, we  consider $\Omega= \sin^2\theta\, \left(\Omega_{\SO(2)}-\Omega_{\SO(n+1)}-r^2\right)$
and write $( \ker \Omega ) _{\widetilde{K}}$ for the space of $\widetilde{K}$-finite vectors in 
$\ker \Omega\subset  \mathcal{I}''_{m,r}$. 
Furthermore, we define two $\widetilde{K}$-representations
$(\mathcal{H}^+)_{\widetilde{K}}$ and $(\mathcal{H}^-)_{\widetilde{K}}$ by
\begin{equation*}
(\mathcal{H}^\pm)_{\widetilde{K}}=\bigoplus\limits_{k\geq 0}\mathbb{C}%
e^{\pm i( k-r) \varphi }\otimes \mathcal{H}_{k}(S^n),
\end{equation*}%
where $\mathcal{H}_{k}(S^n)$ is the space of 
degree $k$ harmonic polynomials on $\mathbb{R}^{n+1}$ restricted injectively to $S^n$. 

\begin{theorem}\label{Thm: K-types}
For $n$ odd, as a $\widetilde{K}$-representation, 
\begin{equation*}
( \ker \Omega ) _{\widetilde{K}}\cong
\begin{cases}
(\mathcal{H}^+)_{\widetilde{K}}\oplus (\mathcal{H}^-)_{\widetilde{K}} & 
\mbox{if\quad $m\equiv n-1 \mod 4$,}\\
\quad\quad\quad\ 0 & \mbox{otherwise},
\end{cases}
\end{equation*}
and for $n$ even,
\begin{equation*}
( \ker \Omega ) _{\widetilde{K}}\cong
\begin{cases}
(\mathcal{H}^+)_{\widetilde{K}} & 
\mbox{if\quad $m\equiv -(n-1) \mod 4$,}\\
(\mathcal{H}^-)_{\widetilde{K}} & 
\mbox{if\quad $m\equiv +(n-1) \mod 4$,}\\
\quad 0 & \mbox{otherwise}.
\end{cases}\qquad\quad
\end{equation*}%
\end{theorem}

\begin{proof}
It follows from the definition  of $\mathcal{I}''_{m ,r}$ that
as  a $\widetilde{K}$-representation,
$$\mathcal{I}''_{m ,r} \isom 
\left\{\phi \in \mathcal{C}^{\infty}(S^{1}\times S^{n})
\mid \phi(w\cdot c)=i^{-m} \phi(c) \quad \forall c\in S^{1}\times S^{n}
\right\}.
$$
The space 
of $\widetilde{K}$-finite vectors in $\mathcal{C}^{\infty}(S^{1}\times S^{n})$ 
decomposes as 
$
\mathcal{C}^{\infty}(S^{1}\times S^{n})_{\widetilde{K}}
\ \isom\ \bigoplus_{p,k}\ 
\CC e^{ip\frac{\varphi}{2}} \otimes \mathcal{H}_{k}(S^n),
$
where the sum is over all $p\in \ZZ$ and $k\in \ZZ_{\geq 0}$.
Noting that on $\CC e^{ip\frac{\varphi}{2}} \otimes \mathcal{H}_{k}(S^n)$ the operator
$\Omega_{\SO(2)}-\Omega_{\SO(n+1)}-r^{2}$ acts  
as the scalar
$$ 
 \left(\frac{p}{2}\right)^{2} - k(k+n-1) -r^{2} =  \left(\frac{p}{2}\right)^{2} -
(k-r)^{2},
$$ 
and $w$ acts as $i^p(-1)^k$,
the result follows immediately.
\end{proof}


\section{Unitarity}
\noindent
The \emph{Klein-Gordon inner product} on the space of smooth solutions
(satisfying appropriate decay conditions) to
the wave equation on $\mathbb{R}^{1,n}$  is defined as%
\begin{equation*}
\left\langle f_{1},f_{2}\right\rangle =i\int_{\mathbb{R}^{n}}( 
\overline{\partial _{t}f_{1}}f_{2}-\overline{f_{1}}\partial _{t}f_{2})
\big|_{t=t_0}\,dx\text{.}
\end{equation*}%
The value of the integral above is independent of the choice of $t_0$. In 
the following, we will always choose $t_0=0$.

\begin{theorem}
For $r=\frac{1-n}{2}$, the Klein-Gordon inner product is well defined and $%
\widetilde{G}$-invariant on $\ker \square\subseteq \mathcal{I}%
_{m ,r}^{\prime }$. Moreover, if $m\equiv \mp(n-1) \mod 4$, it is positive definite on the subspace of $%
\mathcal{I}_{m ,r}^{\prime }$ corresponding to $(\mathcal{H}^{+})_{\widetilde{K}}$ 
and negative definite on the subspace of $\mathcal{I}_{m ,r}^{\prime }$
corresponding to $(\mathcal{H}^{-})_{\widetilde{K}}$.
\end{theorem}

\begin{proof}
If $f_1,f_2\in \mathcal{I}'_{m,r}$, (\ref{compact to noncompact}) implies
$\vert \partial_t f_1(0,x) f_2(0,x)\vert 
\leq C\left( 1+\Vert x\Vert ^{2}\right) ^{-n}$
and hence the Klein-Gordon inner product is well-defined 
on $\mathcal{I}_{m ,r}^{\prime }$.
The fact that the Klein-Gordon inner product is  $\fg$-invariant on
$\ker \square\subseteq \mathcal{I}_{m ,r}^{\prime }$
is  easily checked by integration by parts.

For the last statement of the theorem  
it is useful  to calculate the Klein-Gordan inner product in the compact picture.
If $f_1,f_2\in \mathcal{I}'_{m ,r}$ and $F_1,F_2\in \mathcal{I}''_{m ,r}$ 
are the corresponding functions in the compact picture, then
\begin{equation}\label{Klein-Gordon compact}
\langle f_{1},f_{2}\rangle =
i2^{-n}\int_{[0,\pi ]\times S^{n-1}} 
\left(\,\overline{  \partial_{\varphi}F_{1}}\,F_2 -\overline{ F_{1}}\, \partial_\varphi F_2  \right)
\!\big|_{\varphi=0}\sin^{n-1}\theta\, d\theta d\widehat{x},
\end{equation}
where $d\widehat{x}$ is the spherical measure on $S^{n-1}$.
Now consider $0\not=f\in \mathcal{I}_{m ,r}^{\prime }$ corresponding to a function $F\in(\mathcal{H}^{+})_{\widetilde{K}}$
in the $\widetilde{K}$-type $\CC e^{i( k-r) \varphi }\otimes \mathcal{H}_{k}(S^n)$.
Then  $\partial_\phi F= i(k-r)F$ and since $k-r=k+\frac{n-1}{2}>0$, 
(\ref{Klein-Gordon compact}) implies  $\langle f,f\rangle >0$. 
\end{proof}

\section{An Explicit Orthonormal Basis of Solutions}

\noindent
Let $r=\frac{1-n}{2}$. For $l\in \ZZ_{\geq 0}$ and $p\in \ZZ_{>0}$ of the
form $p=2(l+d-r)$ with $d\in  \ZZ_{\geq 0}$, we define a polynomial
$g_{p, l}(t,x)$ of degree $2d$ by
\begin{equation*}
g_{p, l}(t,x)=\lambda(t,x)^{d}\ 
\widetilde{C}^{l-r}_{d}\left(\frac{1-q(t,x)}{\lambda(t,x)}\right),
\end{equation*}
where 
$q(t,x)=-t^{2}+\Vert x\Vert^{2}$, 
$\lambda(t,x)=\left((1-q(t,x))^{2}+4\Vert x\Vert^{2}\right)^{\frac{1}{2}}$, 
and $\widetilde{C}^{l-r}_{d}(s)$
is the normalized Gegenbauer  polynomial of degree $d$ and parameter 
$l-r$. 
Let $h_{l,j}(x)$ be homogeneous harmonic polynomials on $\RR^{n}$ of degree $l$ such that 
the functions $h_{l,j}|_{S^{n-1}}$ form an orthonormal basis for $\mathcal{L}^{2}(S^{n-1})$. Without loss of generality, we may assume
that the functions $h_{l,j}|_{S^{n-1}}$ are real-valued. (We will need this
assumption in the last section.)

\begin{theorem}
For $r=\frac{1-n}{2}$ and $m= -(n-1) \operatorname{mod} 4$,  the functions 
\begin{equation*}
  f_{p,l,j}(t,x)= \frac{1}{2^{l-r}p^{\frac{1}{2}}}\,\frac{g_{p,l}(t,x)h_{l,j}(x)}{\left(\sqrt{(1-it)^{2}+\Vert x\Vert^{2}}\,\right)^{p}}
\end{equation*}
form an orthonormal basis of the subspace of $\ker\Box\subset \mathcal{I}'_{m,r}$
corresponding to $(\mathcal{H}^{+})_{\widetilde{K}}$.
Similarly, the complex conjugate functions $\overline{f}_{p,l,j}$
form an orthonormal basis of the subspace of 
$\ker\Box\subset \mathcal{I}'_{-m,r}$
corresponding to $(\mathcal{H}^{-})_{\widetilde{K}}$.
\end{theorem}

\begin{proof}
We first work in the compact picture.
Noting that $\mathcal{H}_k(S^n)\isom \bigoplus_{l=0}^k\mathcal{H}_l(S^{n-1})$,
it is clear that the $\widetilde{K}$-type
$\CC e^{i(k-r)\varphi} \otimes \mathcal{H}_{k}(S^n)$ 
is spanned by functions of the form $F(\varphi,\theta,\widehat{x})
=e^{i(k-r)\varphi}g(\cos \theta)\, \sin^l \theta\,h_{l,j}(\widehat{x})$,
where $0\leq l \leq k$ and $g(s)$ is a polynomial
such that $g(\cos \theta)\, \sin^l \theta\,h_{l,j}(\widehat{x})\in \mathcal{H}_k(S^n)$. Using (\ref{spherical laplacian})
and substituting $s=\cos \theta$, 
a straightforward calculation shows that the last condition is equivalent 
to $g(s)$ satisfying the differential equation
\begin{equation*}\label{E:Gegenbauer diff eq}
(1-s^{2}) g''(s)-\left(2(l-r)+1\right)s g'(s)+ d\left(2(l-r)+d\right)g(s) =0,
\end{equation*}
where $d=k-l$. 
This equation is a Gegenbauer differential equation 
and has a unique (up to multiple)
nonzero polynomial solution, namely the Gegenbauer polynomial 
$g(s)=C^{l-r}_d(s)$. 
Writing $p=2(k-r)=2(l+d-r)$, we now define 
\begin{equation*}
F_{p,l,j}(\varphi,\theta,\widehat{x})=
p^{-\frac{1}{2}} e^{i p\frac{\varphi}{2}}\widetilde{C}^{l-r}_{d}(\cos \theta)\, \sin^l \theta\,h_{l,j}(\widehat{x}).
\end{equation*} 
Using (\ref{compact to noncompact}), the functions $F_{p,l,j}$ 
are seen to correspond to the functions $f_{p,l,j}$ in the noncompact picture. 
Finally, using (\ref{Klein-Gordon compact}), it is straightforward to check that the  functions $f_{p,l,j}$ form an orthonormal set and hence an orthonormal basis of the subspace of $\ker\Box\subset \mathcal{I}'_{m,r}$
corresponding to $(\mathcal{H}^{+})_{\widetilde{K}}$.
\end{proof}

\section{Energy and Irreducibility}

\noindent
We define another $\sl(2)$-triple $\{z,n^+,n^-\}$ in $\fg_\CC$ by
\begin{equation*}
   z=i(e^+-e^-), \quad n^+= \textstyle{\frac{1}{2}(h-i(e^{+} + e^{-}))}, \quad
    n^-=\textstyle{\frac{1}{2}(h+i(e^{+} + e^{-}))}.
\end{equation*}
Then $z$ lies in the center of $\fk_\CC$ and $\operatorname{ad}(z)$ acts with eigenvalues
$\{-2,0,+2\}$ on $\fg_\CC$. The corresponding eigenspace decomposition 
$\fg_\CC=\fp^-\oplus \fk_\CC\oplus \fp^+$ is the usual complexified Cartan decomposition.
In the compact picture, the action of  $\{z,n^+,n^-\}$ is given by the formulas
$z =-2i\partial_{\varphi}$ and $n^{\pm} =  e^{\pm i\varphi}\left(r\cos \theta
\pm i  \cos \theta \, \partial_{\varphi}
-  \sin \theta \, \partial_{\theta}\right).$
In particular, this shows that the decomposition of 
$(\ker \Omega)_{\widetilde{K}}\subset \mathcal{I}''_{m,r}$ into $\widetilde{K}$-types
given by Theorem~\ref{Thm: K-types} is the decomposition into eigenspaces of   $z$, which are also referred to as {\it energy levels}.
Using  properties of Gegenbauer polynomials, it is easy to show that
$z\cdot F_{p,l,j}=p  F_{p,l,j}$ and $n^{\pm}\cdot F_{p,l,j} \in \CC F_{p\pm 2,l,j}$ with 
$n^{\pm}\cdot F_{p,l,j}=0$ if and only if $p=\mp2(l-r)$, respectively.
Here for $p<0$,  $F_{p,l,j}=\overline{F}_{-p,l,j}$.
Thus, if for fixed $l\geq 0$ we define 
$V_{l-r}=\operatorname{span}_\CC\{F_{p,l,j} \mid p\geq 2(l-r)\}$ and 
$V^{-(l-r)}=\operatorname{span}_\CC\{F_{p,l,j} \mid p\leq -2(l-r)\}$,
then $V_{l-r}$ is a lowest weight representation of
$\sl(2,\CC)$ with lowest weight $l-r$ and $ V^{-(l-r)}$ is a highest weight representation of $\sl(2,\CC)$ with highest weight
$-(l-r)$. Furthermore, as an $\sl(2,\CC)\times \SO(n)$-representation,
\begin{equation*}
(\mathcal{H}^+)_{\widetilde{K}} \isom \bigoplus_{\!\!\!\!\! l\geq 0\!\!\!\!\!}  \ V_{\ l-r} \otimes \mathcal{H}_l(S^{n-1}) \quad 
\end{equation*}
and
\begin{equation*} 
(\mathcal{H}^-)_{\widetilde{K}} \isom \bigoplus_{\!\!\!\!\! l\geq 0\!\!\!\!\!}  \ V^{-(l-r)} \otimes \mathcal{H}_l(S^{n-1}).
\end{equation*}
From these observations it follows that the $(\fg_\CC, \widetilde{K})$-module
$(\mathcal{H^+})_{\widetilde{K}}$
is an irreducible lowest weight representation with lowest weight vector $e^{-i r\varphi}$
and 
$(\mathcal{H^-})_{\widetilde{K}}$
is an irreducible highest weight representation with lowest weight vector $e^{ir\varphi}$.
To summarize we then have the following result.

\begin{theorem}
Let $r=\frac{1-n}{2}$. Then for $n$ odd, 
$\ker \square\subset \mathcal{I}_{m ,r}^{\prime }$ completes to a 
$\widetilde{G}$-representation%
\begin{equation*}
\isom
\begin{cases}
\mathcal{H}^{+}\oplus \mathcal{H}^{-} & 
\mbox{if\quad $m\equiv n-1 \mod 4$,}\\
\quad\ \ \ 0 & \mbox{otherwise},
\end{cases}
\end{equation*}
and for $n$ even,
\begin{equation*}
\isom
\begin{cases}
\mathcal{H}^{+} & 
\mbox{if\quad $m\equiv -(n-1) \mod 4$,}\\
\mathcal{H}^{-} & 
\mbox{if\quad $m\equiv +(n-1) \mod 4$,}\\
\ 0 & \mbox{otherwise},
\end{cases}\quad
\end{equation*}%
where $\mathcal{H}^{+}$ is a unitary lowest weight representation 
(or positive energy representation) of lowest
weight $-r\omega _{0}$ and $\mathcal{H}^{-}$ is a unitary highest weight
representation (or negative energy representation) of highest weight $r\omega_{0}$, where $\omega_0$ is the fundamental weight corresponding
to the noncompact simple root.
\qed
\end{theorem}

\section{Classical Solutions}

\noindent
It is easy to see that 
every element $f\in \mathcal{H}^+$ is a weak solution to the wave equation.
Here we outline a proof that every classical real-valued solution
is the real part of a unique continuous element in the representation $\mathcal{H}^+$.

\begin{theorem}
Suppose $ \Phi, \Psi \in \mathcal{S}\left( \mathbb{R}^{n}\right) $ are 
real-valued Schwartz functions and let $u\in \mathcal{C}^{\infty}(\RR^{1,n})$ 
be the  solution to the Cauchy problem
\begin{equation*}
\Box u=0, \qquad
u(0,x) =\Phi (x), \qquad \partial_{t}u(0,x) =\Psi (x).
\end{equation*}
Then there is a unique continuous function $f\in \mathcal{H}^{+}$
such that $u=\operatorname{Re}f.$ Explicitly,
$f=\sum_{p,l,j} c_{p,l,j} f_{p,l,j}$ with  $c_{p,l,j}=2\,\langle f_{p,l,j},u\rangle
=2i\!\int_{\mathbb{R}^{n}}( 
\overline{\partial _{t}f_{p,l,j}(0,x)}\Phi(x)-\overline{f_{p,l,j}(0,x)}\Psi(x))dx$.
\end{theorem}

\begin{proof}
Let $r=\frac{1-n}{2}$ and $m=-(n-1) \mod 4$.
The $\fg$-actions on $\mathcal{I}_{m,r}'$ and $\mathcal{I}_{m,r}''$ extend to 
$\fg$-actions on $\mathcal{C}^\infty(\RR^{1,n})$ and 
$\mathcal{C}^\infty(\RR\times\RR\times S^{n-1})$, respectively.
Moreover, the intertwining map $\mathcal{I}_{m,r}'\rightarrow \mathcal{I}_{m,r}''$
given by (\ref{noncompact to compact}) extends to an intertwining map 
$T: \mathcal{C}^\infty(\RR^{1,n}) \rightarrow  
\mathcal{C}^{\infty}\left(\{(\varphi ,\theta ,\widehat{x})\in \RR\times\RR\times S^{n-1}\mid\cos \varphi +\cos
\theta \not=0\}\right)$.
Let $U=T(u)$ and view $U$ as a function on an open dense subset of $\RR \times S^{n}$.
By the assumptions of the theorem, it easily follows that
$({z}^{k}\cdot U)\vert_{\varphi=0}=T({z}^{k}\cdot u)\vert_{\varphi=0}$ is in $\mathcal{L}^2(S^n)$
with values in  $i^k\RR$.
Define $G_{p,l,j}=2^{-r}p^{-\frac{1}{2}} F_{p,l,j}$. Then the real-valued functions
$G_{p,l,j}\vert_{\varphi=0}$ form an orthonormal basis  of 
$\mathcal{L}^2(S^n)$. 
Hence we can write 
$({z}^{k}\cdot U)\vert_{\varphi=0}=\sum_{p,l,j}a_{p,l,j}^{(k)} G_{p,l,j}\vert_{\varphi=0}$
with coefficients $a_{p,l,j}^{(k)}\in i^k\RR$ such that $\sum_{p,l,j} |a_{p,l,j}^{(k)}|^2<\infty$.
A direct calculation in the compact picture shows $c_{p,l,j} =2^{r}p^{-\frac{1}{2}}(
p\,a_{p,l,j}^{(0)}+a_{p,l,j}^{(1)})$.
More generally, $p^k c_{p,l,j}= 2\langle z^k \cdot  f_{p,l,j},u \rangle=2\langle f_{p,l,j}, z^k \cdot u\rangle 
=2^{r}p^{-\frac{1}{2}}(
p\,a_{p,l,j}^{(k)}+a_{p,l,j}^{(k+1)})$
and hence
$a_{p,l,j}^{(k)}=p^{k}a_{p,l,j}^{(0)}$ for $k$ even and 
$a_{p,l,j}^{(k)}=p^{k-1}a_{p,l,j}^{(1)}$ for  $k$ odd.
It follows that 
$\sum_{p,l,j} p^N |a_{p,l,j}^{(0)}|^2<\infty$ and 
$\sum_{p,l,j} p^N |a_{p,l,j}^{(1)}|^2<\infty$ for every $N\geq 0$.
Since $a_{p,l,j}^{(0)}\in \RR$ and $a_{p,l,j}^{(1)}\in i\RR$, 
we have  $\sum_{p,l,j} p^N |c_{p,l,j}|^2<\infty$
for every $N\geq 0$.


To show that $\sum_{p,l,j} c_{p,l,j}f_{p,l,j}$ converges uniformly 
we work  in the compact picture.
Noting that $\left\vert G_{p,l,j}\right\vert $ is
independent of $\varphi $, \cite[Corollary 2.9]{SW}  shows
\begin{equation*}
\sum_{l,j}\left\vert G_{p,l,j}\right\vert ^{2}=
\frac{\dim \mathcal{H}_{\frac{p}{2}+r}(\RR^{n+1})}{\operatorname{Surface\ Area}
(S^{n})} \leq C  p^{n-1},
\end{equation*}
where  $C$ is a constant only depending on $n$.
Thus, by H\"older's inequality,
\begin{equation*}
\sum_{p,l,j}\left\vert c_{p,l,j}\,p^{-\frac{1}{2}}G_{p,l,j}\right\vert 
\leq \bigg( \sum_{p,l,j}p^{n}\left\vert c_{p,l,j}\right\vert ^{2}\bigg)
^{\frac{1}{2}}\bigg( \sum_{p}p^{-n-1}\sum_{l,j}\left\vert
G_{p,l,j}\right\vert ^{2}\bigg) ^{\frac{1}{2}}<\infty.
\end{equation*}
This shows that 
$\sum_{p,l,j} c_{p,l,j}F_{p,l,j}=
2^{-r}\sum_{p,l,j} c_{p,l,j}p^{-\frac{1}{2}} G_{p,l,j}$ converges uniformly
to a continuous function $F$
and hence $\sum_{p.l,j} c_{p,l,j} f_{p,l,j}$ converges uniformly
to a continuous function $f$. Moreover, 
$\Re F\big\vert_{\varphi=0}
=\Re \sum_{p,l,j} (a_{p,l,j}^{(0)}+p^{-1}a_{p,l,j}^{(1)}) G_{p,l,j}\big\vert_{\varphi=0}
=\sum_{p,l,j} a_{p,l,j}^{(0)} G_{p,l,j}\big\vert_{\varphi=0}
=U\vert_{\varphi=0}$ and hence $\Re f \vert_{t=0}= u\vert_{t=0}=\Phi$. 

A subtle argument using a variation of Bernstein's inequality and again
H\"older's inequality shows 
\begin{equation*}
\sum_{p,l,j}\left\vert c_{p,l,j}\,p^{-\frac{1}{2}}\partial_\theta G_{p,l,j}\right\vert 
<\infty.
\end{equation*}
Noting that $e^+=-r\sin \varphi\cos\theta -(1+\cos \varphi\cos \varphi)\partial_\varphi + 
\sin \varphi \sin \theta\, \partial_\theta$ in the compact picture, 
we find that
$\sum_{p,l,j} c_{p,l,j} (e^+ \cdot F_{p,l,j})$ 
converges uniformly. Since $e^+=-\partial_t$ in the noncompact
picture, this implies that 
$\partial_t f = \sum_{p,l,j} c_{p,l,j} \partial_t f_{p,l,j}$
is continuous. Moreover, it follows as above that $\Re \partial_t f \vert_{t=0}= \partial_t u\vert_{t=0}=\Psi$.

Since $\sum_{p,l,j} |c_{p,l,j}|^2<\infty$, $f\in \mathcal{H}^+$.
Since $\Re f $ is a  weak solution to the wave equation with 
$\Re f \vert_{t=0}=\Phi$ and $\Re \partial_t f \vert_{t=0}=\Psi$
it follows that $u=\Re f$.
 \end{proof}

\medskip

\def\germ{\mathfrak}\def\cprime{$'$}\def\scr{\mathcal}

\renewcommand{\baselinestretch}{1.2}

\end{document}